\documentclass[11pt]{amsart}
\usepackage{amssymb,color}
\usepackage{amsmath}
\usepackage{latexsym}
\usepackage{verbatim}
\usepackage[left=1.3in,right=1.3in,top=1in,bottom=1in]{geometry}
\usepackage{color}

\date{\today}

\newcommand{\C}{{\mathbb C}}
\newcommand{\Z}{{\mathbb Z}}

\newcommand{\D}{{\mathbb D}}
\newcommand{\T}{{\mathbb T}}

\newcommand{\N}{{\mathbb N}}

\newtheorem{theorem}{Theorem}[section]
\newtheorem{lemma}[theorem]{Lemma}
\newtheorem{prop}[theorem]{Proposition}

\newtheorem{coro}[theorem]{Corollary}

\theoremstyle{definition}
\newtheorem*{remark}{Remark}
\theoremstyle{definition}

\def\Proof{
\noindent \it Proof.\ \ \rm}
\def\qedbox{$\rlap{$\sqcap$}\sqcup$}

\numberwithin{theorem}{section}
\numberwithin{equation}{section}

\begin{document}

\author[Y.\ Lin]{Yanxue Lin}
\address{School of Mathematical Sciences,  Ocean University of China,
Qingdao
266100, P.R.China}
\email{yanxueLin@aliyun.com}

\author[D.\ Piao]{Daxiong Piao}
\address{School of Mathematical Sciences,  Ocean University of China,
Qingdao
266100, P.R.China}
\email{dxpiao@ouc.edu.cn}

\author[S.\ Guo]{Shuzheng Guo}
\address{School of Mathematical Sciences,  Ocean University of China,
Qingdao
266100, P.R.China}
\email{guoshuzheng@ouc.edu.cn}

\title[Anderson localization for skew-shift Verblunsky coefficients]{Anderson localization for the quasi-periodic CMV matrices with Verblunsky coefficients defined by the skew-shift}

\begin{abstract}
In this paper, we study quasi-periodic CMV matrices with Verblunsky coefficients given by the skew-shift. We prove the positivity of Lyapunov exponents and Anderson localization for most frequencies, which establish the analogous results of one-dimensional Schr\"{o}dinger operators proved by Bourgain, Goldstein and Schlag \cite{Bourgain03}.
\end{abstract}





\maketitle


\section{Introduction}

This paper is concerned with the positivity of the Lyapunov exponents and Anderson localization for quasi-periodic CMV matrices with Verblunsky coefficients given by the skew-shift. Anderson localization is one of the central topics in spectral theory, which describes insulating behavior in the sense that quantum states are essentially localized in a suitable bounded region for all time. In mathematics, Anderson localization often means that the corresponding operator has pure point spectrum with exponentially decaying eigenfunctions. We refer the reader to \cite{LTW09} for the history of the Anderson localization.

Fruitful results are obtained on Anderson localization for Schr\"{o}dinger operators, see for examples \cite{AJ08, Bour07, Bourgain03, Bourgain01, BGS02, Bourgain02, ChuSinai89, CW21, Jito99, Klein05} and references therein. Due to the well-known analogy between the OPUC theory and the theory of discrete one-dimensional Schr\"{o}dinger operators, people expect the CMV matrix analogs of the existing Schr\"{o}dinger operator results. However, the methods on one side do not directly apply for the other side. For CMV matrices, the Anderson localization results are still seldom.

In this paper, we focus our attention to the quasi-periodic CMV matrices. In 2018, Wang and Damanik \cite{Wang01} obtained Anderson localization for the quasi-periodic CMV matrices with Verblunsky coefficients defined by the shift, which established the CMV matrix analog of a result proved by Bourgain and Goldstein \cite{BGS02}. More recently, by following the method of \cite{BGS02}, Cedzich and Werner \cite{CW21} proved Anderson localization for one-dimensional quantum walks placed into homogenous electric fields, which implies Anderson localization for CMV matrices with a particular choice of skew-shift Verblunsky coefficients. Skew-shift model is one of the most important ones in the spectral theory of Sch\"odinger operators and has been widely studied; see, for examples \cite{ABD, Bourgain03, DD, KrJFA, TaoK} and references therein. Inspired by these works, we aim to get the CMV matrix analog of the results proved by Bourgain, Goldstein and Schlag \cite{Bourgain03}.

In the proof of Anderson localization of quasi-periodic Schr\"{o}dinger operators \cite{BGS02,Bour07,Bourgain01,Bourgain03}, some tools such as the avalanche principle and semi-algebraic sets are often used. It is natural to consider whether one can prove Anderson localization for quasi-periodic CMV matrices with skew-shift Verblunsky coefficients by using the analogous tools and methods of \cite{Bourgain03}. This paper provides a positive answer to this question, but the related measure estimates and derivations are much more complicated.

The rest of this paper is structured as follows. We describe the setting and the main results in Section 2. The large deviation estimate and positivity of Lyapunov exponents are shown in Section 3. The main Anderson localization result is then proved in Section 4.  The appendix contains material that is crucial to our work in the main part of the paper.

\medskip

\section{Setting and results}

In this section, we describe the setting in which we work and state the main results. One may consult \cite{S04, S05} for the general background.

We consider the \emph{extended CMV matrix}, which is a pentadiagonal unitary operator on $\ell^2(\Z)$ with a repeating $2\times 4$ block structure of the form
\begin{equation}{\label{1.2}}
\mathcal{E}=\left(
\begin{matrix}
\cdots&\cdots&\cdots&\cdots&\cdots&\cdots&\cdots&\\
\cdots&-\overline{\alpha}_0\alpha_{-1}&\overline{\alpha}_1\rho_{0}&\rho_1\rho_0&0&0&\cdots&\\
\cdots&-\rho_0\alpha_{-1}&-\overline{\alpha}_1\alpha_{0}&-\rho_1\alpha_0&0&0&\cdots&\\
\cdots&0&\overline{\alpha}_2\rho_{1}&-\overline{\alpha}_2\alpha_{1}&\overline{\alpha}_3\rho_2&\rho_3\rho_2&\cdots&\\
\cdots&0&\rho_2\rho_{1}&-\rho_2\alpha_{1}&-\overline{\alpha}_3\alpha_2&-\rho_3\alpha_2&\cdots&\\
\cdots&0&0&0&\overline{\alpha}_4\rho_3&-\overline{\alpha}_4\alpha_3&\cdots&\\
\cdots&\cdots&\cdots&\cdots&\cdots&\cdots&\cdots&
\end{matrix}
\right),
\end{equation}
where $\alpha_n\in\D:=\{z\in\C:|z|<1\}$ is so-called \emph{Verblunsky coefficient} and $\rho_n=\sqrt{1-|\alpha_n|^2}$ for all $n\in\Z$. Setting $\alpha_{-1}=-1$, the matrix decouples into two half-line matrices. The matrix on the right half-line takes the form
\begin{equation}{\label{1.1}}
\mathcal{C}=\left(
\begin{matrix}
\overline{\alpha}_0&\overline{\alpha}_1\rho_{0}&\rho_1\rho_0&0&0&\cdots&\\
\rho_0&-\overline{\alpha}_1\alpha_{0}&-\rho_1\alpha_0&0&0&\cdots&\\
0&\overline{\alpha}_2\rho_{1}&-\overline{\alpha}_2\alpha_{1}&\overline{\alpha}_3\rho_2&\rho_3\rho_2&\cdots&\\
0&\rho_2\rho_{1}&-\rho_2\alpha_{1}&-\overline{\alpha}_3\alpha_2&-\rho_3\alpha_2&\cdots&\\
0&0&0&\overline{\alpha}_4\rho_3&-\overline{\alpha}_4\alpha_3&\cdots&\\
\cdots&\cdots&\cdots&\cdots&\cdots&\cdots&
\end{matrix}
\right),
\end{equation}
and is known as a \emph{standard} or \emph{half-line CMV matrix}.

\medskip

For convenience, let us describe the setting first.
In this paper, we consider a sequence of Verblunsky coefficients generated by an analytic function $\alpha(\cdot,\cdot):\mathbb{T}^2\to\mathbb{D}$, i.e.,
$\alpha_n(x,y)=\lambda \alpha(T_{\omega}^{n} (x,y))$, where $\lambda\in(0,1)$ is a coupling constant, $T_{\omega}(x,y)=(x+y,y+\omega)\in\mathbb{T}^2$ is the skew-shift on the two-dimensional torus $\mathbb{T}^{2}$, and $(x,y)\in\mathbb{T}^2$, $\omega\in\mathbb{T}$ are called \emph{phase} and \emph{frequency} respectively. We assume that the sampling function $\alpha(x,y)$ satisfies
\begin{equation}\label{e.assumptionforalpha}
\int_{\mathbb{T}^{2}}\log(1-|\alpha(x,y)|)d\mu > -\infty.
\end{equation}
The number $\omega$ will be assumed to be Diophantine in the sense that
\begin{equation}\label{e.diophantine}
\|n\omega\|\ge \varepsilon n^{-1} (1+\log n)^{-2} \quad \textrm{for any } n\in \N,
\end{equation}
where $\varepsilon>0$ is some arbitrary but fixed small number. Let $\Omega_{\varepsilon}$ be the set of those $\omega$ that satisfy \eqref{e.diophantine}. It is clear that
\[
\mathrm{mes}[\mathbb{T}\backslash\Omega_{\varepsilon}]<C\varepsilon
\]
with an absolute constant $C$.


According to the analyticity, $\alpha(x,y)$ can be boundedly extended to a strip (see \cite[Chap. 2, Theroem 6]{BM48}) $$\mathcal{D}_{h_1}\times\mathcal{D}_{h_2}:=\{w_1\in\mathbb{C}:1-h_1<|w_1|<1+h_1\}\times\{w_2\in\mathbb{C}:1-h_2<|w_2|<1+h_2\},$$
where $h_1,h_2>0$, with the norm
$$\|\alpha\|_{h_1,h_2}=\sup_{(w_1,w_2)\in\mathcal{D}_{h_1}\times\mathcal{D}_{h_2}}|\alpha(w_1,w_2)|.$$

\medskip

Suppose $\mu$ is a non-trivial (i.e., not finitely supported) probability measure on the unit circle $\partial \mathbb{D} = \{ z \in \mathbb{C} : |z| = 1 \}$, which means the support of $\mu$ contains infinitely many points.
By the non-triviality assumption, the functions $1, z, z^2, \cdots$ are linearly independent in the Hilbert space $\mathcal{H} = L^2(\partial\mathbb{D}, d\mu)$, and hence one can form, by the Gram-Schmidt procedure, the \emph{monic orthogonal polynomials} $\Phi_n(z)$, whose Szeg\H{o} dual is defined by $\Phi_n^{*} = z^n\overline{\Phi_n({1}/{\overline{z}})}$. The Verblunsky coefficients obeys the following equation
\begin{equation}\label{eq01}
\Phi_{n+1}(z) = z \Phi_n(z) - \overline{\alpha}_n \Phi_n^*(z), \qquad \textrm{ for } n\in \Z_+,
\end{equation}
which is known as the \emph{Szeg\H{o} recurrence}.

Consider instead the orthonormal polynomials
$$\varphi(z, n)=\frac{\Phi_n(z)}{\|\Phi_n(z)\|_{\mu}},$$
where $\|\cdot\|_{\mu}$ is the norm of $L^{2}(\partial\mathbb{D},d\mu)$, the Szeg\H{o} recurrence  \eqref{eq01} is equivalent to the following one,
\begin{equation}\label{eq01b}
\rho_n(x,y) \varphi_{n+1}(z)  = z \varphi_{n}(z)  - \overline{\alpha}_n (x,y)\varphi^*_{n}(z), \textrm{ for } n\in \Z_+ ,
\end{equation}
where $\rho_n(x,y) = \sqrt{1 - |\alpha_n(x,y)|^2}$.

The Szeg\H{o} recurrence can be written in a matrix form as follows:
\begin{equation*}
\left(
\begin{matrix}
\varphi_{n+1}\\
\varphi^{*}_{n+1}
\end{matrix}
\right)
=
S^{z}(T^{n}_{\omega}(x,y))
\left(
\begin{matrix}
\varphi_{n}\\
\varphi^{*}_{n}
\end{matrix}
\right), \textrm{ for } n\in\Z_+,
\end{equation*}
where
\begin{equation*}
S^{z}(x,y)
=
\frac{1}{\rho_n(x,y)}
\left(
\begin{matrix}
z & -\overline{\alpha}_n(x,y)\\
-\alpha_n(x,y) z & 1
\end{matrix}
\right).
\end{equation*}

Since det$S^{z}(x,y)=z$, we study the determinant 1 matrix instead
\begin{equation*}
M^{z}(x,y)
=
\frac{1}{\rho(x,y)}
\left(
\begin{matrix}
\sqrt{z} &\frac{ -\overline{\alpha}(x,y)}{\sqrt{z}}\\
-\alpha(x,y)\sqrt{z} & \frac{1}{\sqrt{z}}
\end{matrix}
\right),
\end{equation*}
which is called the  \emph{Szeg\H{o} cocycle map}.
Then the \emph{$n$-step transfer matrix} is defined by
\begin{equation}\label{015}
M^{z}_n(x,y)   = \prod^{0}_{j=n-1}M^{z}(T^{j}_{\omega}(x,y)).
\end{equation}
It follows from \eqref{015} that
$$M^{z}_{n_{1}+n_{2}}(x,y)   =M^{z}_{n_{2}}(T^{n_{1}}_{\omega}(x,y))M^{z}_{n_{1}}(x,y), $$
and hence
\begin{equation}\label{016}
\log\|M^{z}_{n_{1}+n_{2}}(x,y)\|  \leq \log\|M^{z}_{n_{1}}(x,y)\|+ \log\|M^{z}_{n_{2}}(T^{n_{1}}_{\omega}(x,y))\|.
\end{equation}

Integrating \eqref{016} over $\T^2$, we get
$$L_{n_{1}+n_{2}}(z)\leq\frac{n_{1}}{n_{1}+n_{2}}L_{n_{1}}(z)+\frac{n_{2}}{n_{1}+n_{2}}L_{n_{2}}(z),$$
where
$$L_{n}(z)=\frac{1}{n}\int_{\mathbb{T}^{2}}\log\|M^{z}_n(x,y) \|dxdy.$$
This implies that
$$L_{n}(z)\leq L_{m}(z) \quad if \quad  m < n, m \mid n,$$
and
$$L_{n}(z)\leq L_{m}(z)+C\frac{m}{n} \quad if \quad   m < n. $$

Notice that \eqref{016} means that $\log\|M^{z}_n(x,y) \|$ is subadditive. Therefore, the \emph{Lyapunov exponent} $L(z)$ defined by
$$L(z)=\lim_{n\rightarrow \infty}L_{n}(z)$$
exists.

\medskip

With $\{\alpha_n\}_{n\in\Z}$ and $T_{\omega}(x,y)$ as above, a CMV matrix $\mathcal{E}_{\omega}(x,y)$ can be generated dynamically. Our first main result concerns the positivity of Lyapunov exponents.
\begin{theorem}{\label{mainth1}}
Fix $\varepsilon>0$ small and let $\omega\in \Omega_\varepsilon$, see \eqref{e.diophantine}. Let $\alpha_n(x,y)=\lambda \alpha(T_{\omega}^{n} (x,y))$, where $\alpha(x,y)$ is a nonconstant trigonometric polynomial on $\mathbb{T}^2$. For all $z\in \partial\D$, there exists a $\lambda_0\in(0,1)$ such that for every $\lambda\in(\lambda_{0},1)$ and all $(\omega,x,y)\in \mathbb{T}^3$ up to a set of measure $\varepsilon$,
\[
L(z)\ge -\frac{1}{4} \log(1-\lambda^2).
\]
\end{theorem}

\medskip

Based on this theorem, we can obtain the following Anderson localization result.
\begin{theorem}{\label{mainth2}}
Let $\lambda, \omega, x, y$ and $z$ as in Theorem \ref{mainth1}, the CMV operator $\mathcal{E}_{\omega}(x,y)$ displays Anderson localization.                                                                                                                                                                                                                                                                                                                                                                                                                                                                                                                                                                                                                                                                                                                                                                                                                                                                                                                                                                                                                                                                                                                                                                                                                                                                                                                                                                                                                                                                                                                                                                                                                                                                                                                                                                                                                                                                                                                                                                                                                                                                                                                                        \end{theorem}

\begin{remark}

\begin{enumerate}
  \item It is well-known that $\mathcal{C}_{[a,b]}=\mathcal{E}_{[a,b]}$ whenever $1\le a\le b$ and $\mathcal{C}_{[0,b]}=\mathcal{E}_{[0,b]}$ with $\alpha_{-1}=-1$. By running the similar proof, one can easily obtain the analogous results for the half-line case.
  \item From \cite[Theorem 3.4] {Damanik01}, we see that if every generalized eigenfunction $\xi=\{\xi_n\}_{n\in \mathbb{Z}}$ of $\mathcal{E}_\omega$ decays exponentially, then the operator $\mathcal{E}_\omega(x,y)$ displays Anderson localization. That is, it suffices to show that if $z\in \partial\D$ and $\xi=\{\xi_n\}_{n\in \mathbb{Z}}$ satisfying
$$|\xi_n|\lesssim|n|^C$$
and
$$\mathcal{E}_\omega \xi =z\xi,$$
then
$$|\xi_n|\lesssim e^{-c|n|},~~{\rm for ~ some}~ c>0.$$
Where $\lesssim$ denotes inequality up to a multiplicative constant.
\end{enumerate}
\end{remark}




\medskip

\section{A large deviation theorem and positive Lyapunov exponents}

Since $M^z(x,y)$ is conjugate to a $\mathrm{SL}(2,\mathbb{R})$ matrix $A^{z}(x,y)$, that is,
$A^z(x,y)=Q^*M^z(x,y)Q\in\mathrm{SL}(2,\mathbb{R})$, where
$$Q=-\frac{1}{1+i}
\left(
\begin{matrix}
1&-i\\
1&i
\end{matrix}
\right)\in \mathbb{U}(2).$$
A direct calculation implies $\log{\|M^z(x,y)\|}$ is subharmonic, so $\log{\|A^z(x,y)\|}$ is also subharmonic.
By conjugacy, $\|A^z(x,y)\|=\|M^z(x,y)\|$ and $\|A^z_n(x,y)\|=\|M^z_n(x,y)\|$, where
\begin{equation}{\label{Wang2.4}}
A_n^z(x,y)=\prod_{j=n-1}^0 A^z(T_{\omega}^j(x,y)).
\end{equation}

As usual,
\begin{equation*}
L_{n}(z)=\frac{1}{n}\int_{\mathbb{T}^{2}}\log\|A^{z}_n(x,y) \|dxdy,
\end{equation*}
and $L(z)=\lim_{n\rightarrow \infty}L_{n}(z)$ denotes the Lyapunov exponent. Introduce a scaling factor
\begin{equation}\label{2.61}
P(z)=\log\left(\sup_{j\in\Z}\frac{1}{\rho_j}+C_{\alpha}+(1-\lambda^2)^{-1}+|z|\right)\geq 1,
\end{equation}
where $C_{\alpha}$ is a constant depending only on $\alpha$ so that for all $n$
\begin{equation}\label{e.subharmonic}
\sup_{w_1\in\mathcal{D}_{h_1}}\sup_{w_2\in\mathcal{D}_{h_2}} \frac{1}{n} \log \|M_n(w_1,w_2;z)\|\le P(z).
\end{equation}

\subsection{The main inductive step}
The following lemma provides the inductive step in the proof of the large deviation theorem. 
\begin{lemma}\label{le.BGSle2.8}
 Fix $\varepsilon>0$ small and let $\omega\in\Omega_{\varepsilon}$, see \eqref{e.diophantine}. Suppose $n$ and $N>n$ are positive integers such that
\begin{equation}\label{e.LDTforn}
\mathrm{mes}\left[(x,y)\in\mathbb{T}^2: \left|\frac{1}{n}\log\|M_n(x,y;z)\|-L_n(z)\right|>\frac{\gamma}{10}P(z)\right]\le N^{-10},
\end{equation}
and
\begin{equation}\label{e.LDTfor2n}
\mathrm{mes}\left[(x,y)\in\mathbb{T}^2: \left|\frac{1}{2n}\log\|M_n(x,y;z)\|-L_{2n}(z)\right|>\frac{\gamma}{10}P(z)\right]\le N^{-10},
\end{equation}
where $\gamma >0$ small enough.
Assume that
\begin{equation}\label{e.minofLnL2n}
\min(L_n(z),L_{2n}(z))\ge \gamma P(z),
\end{equation}
\begin{equation}\label{e.diffofLnL2n}
L_n(z)-L_{2n}(z)\le \frac{\gamma}{40}P(z),
\end{equation}
\begin{equation}\label{e.Nleqn2}
9\gamma n P(z) \ge 10\log(2N) \textrm{ and } n^2\le N.
\end{equation}
Then
\begin{equation}\label{e.ineforLNanddiff}
\begin{split}
& L_N \ge \gamma P(z)-2(L_{n}-L_{2n})-C_0nN^{-1}P(z),\\
& L_N-L_{2N}  \le C_0nN^{-1}P(z).
\end{split}
\end{equation}
Moreover, for any $\sigma< \frac{1}{24}$, there is a $\tau=\tau(\sigma)>0$ so that
\begin{equation}\label{e.LDTforN}
\begin{split}
\mathrm{mes}\left[(x,y)\in\mathbb{T}^2: \left|\frac{1}{N}\log\|M_N(x,y;z)\|-L_N(z)\right| >P(z)N^{-\tau}\right]\\
& \quad \le C \exp(-N^{\sigma})
\end{split}
\end{equation}
with some constant $C=C(N,\sigma,\varepsilon)$.
\end{lemma}
\Proof
Denote the set on the left-hand side of \eqref{e.LDTforn} by $\mathcal{B}_{n}$ and the set on the left-hand side of \eqref{e.LDTfor2n} by $\mathcal{B}_{2n}$. For any $(x,y)\in \mathbb{T}^2\backslash \mathcal{B}_{n}$,

\begin{equation*}
\|M_{n}(x,y;z)\|\geq \exp(\gamma n P-\frac{\gamma}{10}nP)=\exp(\frac{9\gamma}{10}n P):=\mu.
\end{equation*}
By \eqref{e.Nleqn2}, $\mu \geq 2N$, which corresponds to the first condition \eqref{2.8} of the avalanche principle (see Lemma \ref{le.AP}). Furthermore, for any $(x,y) \not\in \mathcal{B}_{n}\bigcup T^{-n}\mathcal{B}_{n}\bigcup \mathcal{B}_{2n} $,\eqref{e.LDTforn}-\eqref{e.diffofLnL2n} imply
\begin{equation}\label{e.avalanche02}
\begin{split}
&\log\|M_{n}\circ T^{n}(x,y)\| \left.\right.+\log\|M_{n}(x,y)\|-\log\|M_{2n}(x,y)\| \\
& \leq 2n(L_{n}-L_{2n})+\frac{4\gamma}{10}n P(z)\leq \frac{9\gamma}{20}n P(z) =\frac{1}{2}\log\mu
\end{split}
\end{equation}
and that is the second condition \eqref{2.9} of the avalanche principle. In fact, since
\begin{equation*}
\left|\frac{1}{n}\log\|M_n(x,y;z)\|-L_n(z)\right|\leq\frac{\gamma}{10} P(z),
\end{equation*}
\begin{equation*}
\left|\frac{1}{n}\log\|M_n\circ T^n(x,y;z)\|-L_n(z)\right|\leq\frac{\gamma}{10}P(z),
\end{equation*}
and
\begin{equation*}
\left|\frac{1}{2n}\log\|M_{2n}(x,y;z)\|-L_{2n}(z)\right|\leq\frac{\gamma}{10} P(z),
\end{equation*}
we can get
\begin{equation*}
\frac{1}{n}\log\|M_n(x,y;z)\| \leq \frac{\gamma}{10} P(z)+L_n(z),
\end{equation*}
\[
\frac{1}{n}\log\|M_n\circ T^n(x,y;z)\| \leq \frac{\gamma}{10} P(z)+L_n(z),
\]
and
\begin{equation*}
\frac{1}{2n}\log\|M_{2n}(x,y;z)\| \geq L_{2n}(z)-\frac{\gamma}{10} P(z).
\end{equation*}
It follows that \eqref{e.avalanche02} holds.


\medskip

Applying Lemma \ref{le.AP} $N$ times yields a set $\mathcal{B}_{1}\subset \mathbb{T}^2$ with measure
\begin{equation}\label{e.2.73}
\mathrm{mes}(\mathcal{B}_{1})\leq 4\cdot N\cdot N^{-10}=4N^{-9},
\end{equation}
so that for any $(x,y)\in \mathbb{T}^2\backslash \mathcal{B}_{1}$,
\begin{align}
&\left| \frac{1}{N}\log\|M_{N}(x,y)\| +\frac{1}{N}\sum^{N}_{j=1}\frac{1}{n}\log\|M_{n}\circ T^{j}(x,y)\|  -\frac{2}{N}\sum^{N}_{j=1}\frac{1}{2n}\log\|M_{2n}\circ T^{j}(x,y)\| \right| \notag \\
& \leq \frac{1}{N}(CP(z)n+C\frac{N}{\mu}) \notag \\
& \leq CP(z)nN^{-1}.\label{e.LDT2.74}
\end{align}
Integrating \eqref{e.LDT2.74} over $\mathbb{T}^2$ (seen as $\mathcal{B}_{1} \bigcup \left(\mathbb{T}^2\backslash \mathcal{B}_{1}\right)$) yields
\begin{equation}\label{e.LDT2.75}
|L_{N}+L_{n}-2L_{2n}|\leq CP(z)nN^{-1}+16P(z)N^{-9}.
\end{equation}
Then
\begin{equation*}
\begin{split}
L_{N}&\geq L_n-2(L_{n}-L_{2n})-CP(z)nN^{-1}-16P(z)N^{-9}\\
&\geq \gamma P(z)-2(L_{n}-L_{2n})-C_{0}P(z)nN^{-1}.
\end{split}
\end{equation*}
It follows that the first inequality in \eqref{e.ineforLNanddiff} holds.
Running this procedure with $M_{2N}$ instead of $M_N$, we can get the analogous inequality involving $L_{2N}$. Then one can get the second inequality in \eqref{e.ineforLNanddiff}.


\medskip

Denote
\begin{equation*}
u_{N}(x,y)=\frac{1}{N}\log\|M_{N}(x,y)\|,
\end{equation*}
and similarly with $n$ and $2n$. In view of \eqref{e.subharmonic}, both $u_{n}$ and $u_{2n}$ extend to separately subharmonic functions in both variables such that
\begin{equation*}
\sup_{w_1\in\mathcal{D}_{h_1}}\sup_{w_2\in\mathcal{D}_{h_2}} [|u_{n}(w_1,w_2)|+|u_{2n}(w_1,w_2)|]\leq CP(z).
\end{equation*}

By \cite[p.601]{Bourgain03}, applying the properties of subharmonic functions  and combining Lemma \ref{le.subharmonic01} and Lemma \ref{le.subharmonic02}, we can complete the proof. Indeed, in this case the quantity $B$ from Lemma \ref{le.subharmonic01} satisfies
\begin{equation*}
B \leq C_{\delta}N^{-\frac{1}{12}+2\delta}\log(N^{10})+C_{\delta}^{\frac{1}{4}}N^{\frac{3}{2}-\frac{9}{4}}
\end{equation*}
for the fixed number $\delta>0$ small, which gives us the value of $\sigma$ stated above.
\qedbox

\medskip

\subsection{The initial condition}

For $z\in \partial\mathbb{D}$, define $\hat\varphi_{[a,b]}=\det(z-\mathcal{E}_{[a,b]})$.
From \cite{Wang02}, one can write $M_{n}(x,y;z)$ as follows
\begin{equation}\label{e.detversionforTM}
M_n(x,y;z)=(\sqrt{z})^{-n}(\prod_{j=0}^{n-1}\frac{1}{\rho_j})
\left(
\begin{matrix}
z\hat\varphi_{[1,n-1]} & \frac{z\hat\varphi_{[1,n-1]}-\hat\varphi_{[0,n-1]}}{\alpha_{-1}}\\
z(\frac{z\hat\varphi_{[1,n-1]}-\hat\varphi_{[0,n-1]}}{\alpha_{-1}})^* &(\hat\varphi_{[1,n-1]})^{*}
\end{matrix}
\right).
\end{equation}
\begin{lemma}\label{le.BGSle2.10}
For any positive integer $n$,
\begin{equation}\label{2.83}
\sup_{z}\mathrm{mes}\left[(x,y)\in\mathbb{T}^2: \left|\frac{1}{n}\log\|M_n(x,y;z)\|-L_n(z)\right|>\frac{1}{20}P(z)\right]\le 4n^{-100},
\end{equation}
provided $\lambda$ is sufficiently close to 1. For those $\lambda$ and all $z \in \partial\mathbb{D}$,
\begin{equation*}
L_{n}(z)\geq\frac{1}{2}P(z),~~~~~~~L_{n}(z)-L_{2n}(z)\leq\frac{1}{80}P(z).
\end{equation*}
\end{lemma}
\Proof
Let the matrix $z-\mathcal{E}_{[a,b]}=D_{[a,b]}+B_{[a,b]}$, where $a\le b \in \Z$, and
$$D_{[a,b]}=\textrm{diag}(z+\overline{\alpha}_{a}\alpha_{a-1},z+\overline{\alpha}_{a+1}\alpha_{a},z+\overline{\alpha}_{a+2}\alpha_{a+1},\cdots,z+\overline{\alpha}_{b}\alpha_{b-1}).$$
Without loss of generality, we set $a=0,b=n-1$, then
$$D_{[0,n-1]}=\textrm{diag}(z+\overline{\alpha}_{0}\alpha_{-1},z+\overline{\alpha}_{1}\alpha_{0},z+\overline{\alpha}_{2}\alpha_{1},\cdots,z+\overline{\alpha}_{n-1}\alpha_{n-2}).$$
Clearly, due to $|\alpha_{j}|,|\rho_{j}|<1$, we have $\|B_{[0,n-1]}\|<3$.

\medskip

It is a well-known property of complex-valued real analytic functions $v$ that there exist constants $b>0$ and $C$ depending on $v$ such that
\begin{equation}\label{2.85}
\mathrm{mes}[(x,y)\in\mathbb{T}^2:|v(x,y)-h|<t]\leq Ct^{b},
\end{equation}
for all $-2\|v\|_{\infty}\leq h\leq 2\|v\|_{\infty}$ and $t>0$, see for example \cite[Lemma 11.4]{Goldstein01}.
Therefore,  by~\eqref{2.85}, if $$2\lambda^{2}\sup_{j\in\Z,(x,y)\in\mathbb{T}^2}\left|\overline{\alpha}(T^j(x,y)) \alpha (T ^{j-1}(x,y))\right|:=2\lambda^{2}\|\alpha\|_*\geq 1,$$
one can obtain that
\[
\mathrm{mes}[(x,y)\in\mathbb{T}^2:|\lambda^2\overline{\alpha}(T^j(x,y)) \alpha (T ^{j-1}(x,y))-1 | < e^{-\rho}]<C e^{-b\rho},
\]
where $\rho>0$ is an absolutely constant.

Since
$$\|D_{[0,n-1]}(x,y;\lambda,z)^{-1}\|\leq (1-\lambda^2)\max_{0\leq j\leq n-1}| \lambda^2 (1-\lambda^2) \overline{\alpha}(T^j(x,y))\alpha( T^{j-1}(x,y))+z(1-\lambda^2)|^{-1},$$
\eqref{2.85} implies that

\begin{align}\label{2.88}
&\mathrm{mes}\left[(x,y)\in\mathbb{T}^2:\|D_{[0,n-1]}(x,y;\lambda,z)^{-1}\|> \frac{2}{3}\right]\notag\\
&\leq n~ \mathrm{mes}\left[(x,y)\in\mathbb{T}^2:\left|\lambda^2 (1-\lambda^2)\overline{\alpha}(T(x,y))\alpha(x,y) + z (1-\lambda^2)\right|<\frac{3}{2}(1-\lambda^2)\right]\notag\\
&= n~ \mathrm{mes}\left[(x,y)\in\mathbb{T}^2:\left|z^{-1}\lambda^2(1-\lambda^2)\overline{\alpha}(T(x,y))\alpha(x,y) - (1-\lambda^2)\right|<\frac{3}{2}(1-\lambda^2)\right] \notag \\
&\leq  C n (1-\lambda^2)^b.
\end{align}
Hence,
\begin{equation}\label{2.89}
\mathrm{mes}\left[(x,y)\in\mathbb{T}^2:\|D_{[0,n-1]}(x,y;\lambda;z)^{-1}B_{[0,n-1]}\|\geq 2\right]\leq  C n (1-\lambda^2)^b.
\end{equation}

\medskip

Let $\hat\varphi_{[0,n-1]}=\det(z-\mathcal{E}_{[0,n-1]})$. In view of \eqref{2.88} and \eqref{2.89}, we have
\begin{align}
&\left|\frac{1}{n}\log |\hat\varphi_{[0,n-1]}(x,y;\lambda,z)|\right|\notag\\
&=\left|\frac{1}{n}\log |\det D_{[0,n-1]}(x,y;\lambda,z)| + \frac{1}{n}\log |\det (I+ D^{-1}_{[0,n-1]}(x,y;\lambda,z)B_{[0,n-1]})|\right|\notag\\
&\leq \log(1-\|\alpha\|_{\ast})^{-1}+\log 3, \label{2.90}
\end{align}
up to a set of measure not exceeding
\begin{equation}\label{2.91}
C n (1-\lambda^2)^b.
\end{equation}

\medskip
Now assume that $(3(1-\|\alpha\|_{\ast})^{-1})^{\frac{200}{199}}\leq(1-\lambda^2)^{-1}$, then the right-hand side of \eqref{2.90} is no larger than $\frac{199}{200}\log(1-\lambda^2)^{-1}$.   Choosing
$$1-\lambda^2<n^{-B},$$
for some $B>0$ implies
\begin{align*}
\sup_{2\lambda^{2}\|\alpha\|_*\geq 1}\mathrm{mes} \Big[(x,y)\in\mathbb{T}^2: & \left|\frac{1}{n}\log|\hat\varphi_{[0,n-1]}(x,y;\lambda,z)|-\log (1-\lambda^{2})^{-1}\right|\\
&\geq \frac{1}{200}\log (1-\lambda^2)^{-1}\Big]
\le n^{-100}.
\end{align*}
Running the foregoing steps again, we can also get
\begin{align*}
\sup_{2\lambda^{2}\|\alpha\|_*\geq 1}\mathrm{mes}\Big[(x,y)\in\mathbb{T}^2: & \left|\frac{1}{n}\log|\hat\varphi_{[0,n-1]}-\hat\varphi_{[1,n-1]}|-\log (1-\lambda^{2})^{-1}\right|\\
&\geq  \frac{1}{200}\log (1-\lambda^2)^{-1}\Big]
\le n^{-100}.
\end{align*}
Hence,
\begin{align}
\sup_{2\lambda^{2}\|\alpha\|_*\geq 1}\mathrm{mes}\Big[(x,y)\in\mathbb{T}^2: & \left|\frac{1}{n}\log\|M_n(x,y;z)\|-\log(1-\lambda^2)^{-1}\right|\notag \\
&\geq\frac{1}{199}\log (1-|\lambda|^2)^{-1}\Big]
\le 4n^{-100}. \label{2.92}
\end{align}

In particular,
\begin{equation}\label{2.93}
|L_{n}(z)-\log(1-|\lambda|^2)^{-1}|\leq \frac{1}{199}\log (1-|\lambda|^2)^{-1}+4P(z)n^{-100}\leq\frac{1}{198}P(z).
\end{equation}
Since
\begin{equation*}
\log(1-|\lambda|^2)^{-1}\geq\frac{99}{100}\sup_{2\lambda^{2}\|\alpha\|_*\geq 1}P(z),
\end{equation*}
for $\lambda\rightarrow1$, \eqref{2.93} implies the second statement of the lemma. Hence, we have
\begin{equation*}
\sup_{2\lambda^{2}\|\alpha\|_*\geq 1}\mathrm{mes}\left[(x,y)\in\mathbb{T}^2: \left|\frac{1}{n}\log\|M_n(x,y;z)\|-L_{n}(z)\right|\geq\frac{1}{90}P(z)\right]\le 4n^{-100}.
\end{equation*}

\medskip

If $ 2\lambda^{2}\|\alpha\|_*<1$ and $\lambda \rightarrow 1$, then the set in \eqref{2.83} is empty. In fact, from the definition of $D_{[0,n-1]}$, it follows
\begin{equation*}
\begin{split}
\left|\frac{1}{n}\log|\det D_{[0,n-1]}|\right|&= \left|\frac{1}{n}\sum^{n-1}_{j=0}\log|\lambda^2\overline{\alpha}(T^{j}(x,y))\alpha(T^{j-1}(x,y))+z|\right|.
\end{split}
\end{equation*}
We consider two cases.\\
\textbf{Case 1.} If $|\lambda^2\overline{\alpha}(T^{j}(x,y))\alpha(T^{j-1}(x,y))+z|<1$, then we have
\begin{align*}
|\lambda^2\overline{\alpha}(T^{j}(x,y))\alpha(T^{j-1}(x,y))+z|
\ge & |z|-|\lambda^2\overline{\alpha}(T^{j}(x,y))\alpha(T^{j-1}(x,y))|\\
\ge & |z| -\lambda^2 \|\alpha\|_* \ge \frac{1}{2}.
\end{align*}
\textbf{Case 2.} If $|\lambda^2\overline{\alpha}(T^{j}(x,y))\alpha(T^{j-1}(x,y))+z|>1$, then we have
\begin{align*}
|\lambda^2\overline{\alpha}(T^{j}(x,y))\alpha(T^{j-1}(x,y))+z|
\le  |\lambda^2\overline{\alpha}(T^{j}(x,y))\alpha(T^{j-1}(x,y))| + |z|
\le  \frac{3}{2} <2.
\end{align*}
Thus, we have
\begin{equation}
\left|\frac{1}{n}\log|\det D_{[0,n-1]}(x,y,\lambda,z)|\right|<\log 2< 2,
\end{equation}
and thus
\begin{equation*}
\begin{split}
&\left|\frac{1}{n}\log|\widehat{\varphi}_{[0,n-1]}(x,y,\lambda,z)|\right|\\
&\leq \left|\frac{1}{n}\log|\det D_{[0,n-1]}(x,y,\lambda,z)|\right|+\left|\frac{1}{n} \log |\det(I+D_{[0,n-1]}(x,y;\lambda,z)^{-1}B_{[0,n-1]})|\right|\\
&<2+\log 7\\
&<4
\end{split}
\end{equation*}
which implies that for $\lambda \rightarrow 1$,
\begin{equation*}
\frac{1}{n}\log\|M_n(x,y;z)\|<8+\sup_{j\in\Z}\frac{1}{\rho_j}+C_{\alpha_{-1}}<\frac{1}{200}P(z).
\end{equation*}

Hence
\begin{equation*}
L_{n}(\lambda,z) <\frac{1}{200}P(z),
\end{equation*}
and
\begin{equation*}
\left|\frac{1}{n}\log\|M_n(x,y;z)\|-L_{n}(\lambda,z)\right|<\frac{1}{100}P(z).
\end{equation*}

From what has been discussed above, the lemma follows.
\qedbox

\medskip

\subsection{The proof of the large deviation estimate and Theorem \ref{mainth1}}

Now, we could prove the large deviation theorem and positive Lyapunov exponents.
\begin{prop}
Fix $ \varepsilon > 0 $ small and let $\omega \in \Omega_\varepsilon$, see \eqref{e.diophantine}. Assume $\alpha(x,y)$ is a nonconstant analytic function on $\T^{2}$. Then for all $\sigma<\frac{1}{24}$, there exist $\tau=\tau(\sigma)>0$ and constants $\lambda_1$ and $n_0$ depending only on $\varepsilon$, $v$ and $\sigma$ such that
\begin{equation}
\begin{split}
&\mathrm{mes}\left[(x,y)\in\mathbb{T}^2: \left|\frac{1}{n}\log\|M_n(x,y;z)\|-L_n(z)\right.\right.\left.\left.\right|>P(z)n^{-\tau}\right]\\
& \quad \le C \exp(-n^{\sigma}).
\end{split}
\end{equation}
Furthermore, for those $\lambda$ and all $z\in \partial\mathbb{D}$,
\begin{equation}
L(z)=\lim_{n\rightarrow \infty}L_{n}(z)\geq \frac{1}{4}\log(1-\lambda^2)^{-1}
\end{equation}

\end{prop}
\Proof
Let $\tau=\tau(\sigma)>0$ be as in \eqref{e.LDTforN}, $\lambda \geq \lambda_{0}\vee (1-n_{0}^{-B}):=\lambda_{1}$. Fix and require $n_{0}$ to be sufficiently large at various places.

In view of Lemma \ref{le.BGSle2.10}, the hypotheses of Lemma \ref{le.BGSle2.8} are satisfied with $\gamma=\gamma_{0}=\frac{1}{2}$,
\begin{equation}\label{2.96}
n_{0}^{2} \leq N \leq n_{0}^{5},
\end{equation}
provided
\begin{equation}\label{2.97}
9n_{0} \geq 20\log(2n_{0}^{10})
\end{equation}
cf. \eqref{e.Nleqn2} ($P(z) \geq 1$). It is clear that \eqref{2.97} holds if $n_0$ is large. By Lemma \ref{le.BGSle2.8}, we can get



\begin{equation}\label{2.98}
\begin{split}
&L_{N}\geq (\frac{1}{2}-\frac{1}{40})P(z)-C_{0}P(z)N^{-1}n_{0} \geq \gamma_{1}P(z), \\
&L_{N}-L_{2N} \leq C_{0}P(z)N^{-1}n_{0} \leq \frac{\gamma_1}{40}P(z)
\end{split}
\end{equation}
with $\gamma_{1}=\frac{1}{3}$.

Moreover, with some constant $ C_{1}\geq 1$ depending on $\varepsilon$,
\begin{equation}\label{2.99}
\begin{split}
&\mathrm{mes}\left[(x,y)\in\mathbb{T}^2: \left|\frac{1}{N}\log\|M_N(x,y;z)\|-L_N(z)\right|>P(z)N^{-\tau}\right]\\
& \quad \le C_1 \exp(-N^{\sigma})
\end{split}
\end{equation}
for all $N$ in the range given by \eqref{2.96}.

In particular, \eqref{2.99} implies that
\begin{equation*}
\begin{split}
&\mathrm{mes}\left[(x,y)\in\mathbb{T}^2: \left|\frac{1}{N}\log\|M_N(x,y;z)\|-L_N(z)\right|>P(z)\frac{\gamma _1}{10}\right]\\
&\quad \le C_1 \exp(-N^{\sigma}) \leq \overline{N}^{-10},
\end{split}
\end{equation*}
provided $n_0$ is large and
\begin{equation*}
N^{2}\leq  \overline{N} \leq C_{1}^{-\frac{1}{10}} \exp(\frac{1}{10}N^\sigma).
\end{equation*}

The first inequality was added to satisfy \eqref{e.Nleqn2}. In the view of \eqref{2.96}, one thus has the range
\begin{equation}\label{2.100}
n_{0}^{4} \leq  \overline{N} \leq  \exp(\frac{1}{10}n_{0}^{5\sigma})
\end{equation}
of admissible $ \overline{N}$.
Moreover,
\begin{equation}\label{2.101}
\begin{split}
&L_{\overline{N}} \geq \gamma_{1}P(z)-2C_{0}P(z)N^{-1}n_{0}-C_{0}P(z) \overline{N}^{-1}N,  \\
&L_{\overline{N}}-L_{2\overline{N}} \leq C_{0}P(z)\overline{N}^{-1}N.
\end{split}
\end{equation}

At the next stage of this procedure, observe that the left end-point of the range of admissible indices starts at $n_{0}^8$, which is less than the right end-point of the range \eqref{2.100} (for $n_0$ large). Therefore, from this point on the ranges will overlap and cover all large integers.
To ensure that the process does not terminate, simply note the rapid convergence of the series given by \eqref{2.101}.
\qedbox

\medskip

\section{Localization}

In this section, we present some known results that will be used in the proofs, develop the version of some lemmas we need, and prove the Anderson localization.

\subsection{Green's function estimates}
Let $\{\alpha_n\}_{n\in\Z_+}$ be the sequence of Verblunsky coefficients of a half-line CMV matrix $\mathcal{C}$. Define the unitary matrices
\begin{equation*}
\Theta_{n}=\left[
\begin{matrix}
\overline{\alpha_n} &\rho_n\\
\rho_n & -\alpha_n
\end{matrix}
\right].
\end{equation*}
Denote $\mathcal{L}_{+}$, $\mathcal{M}_{+}$ by
\begin{equation*}
\mathcal{L}_{+}=\left[
\begin{matrix}
\Theta_0 &~ & ~\\
~& \Theta_2 & ~\\
~ & ~& \ddots
\end{matrix}
\right],
\mathcal{M}_{+}=\left[
\begin{matrix}
\mathbf{1} &~ & ~\\
~& \Theta_1 & ~\\
~ & ~& \ddots
\end{matrix}
\right],
\end{equation*}
where $\mathbf{1}$ represents the $1\times1$ identity matrix. It is well known that $\mathcal{C}=\mathcal{L}_{+}\mathcal{M}_{+}$.
Now let $\{\alpha_n\}_{n\in\Z}$ be a bi-infinite sequence of Verblunsky coefficients and denote the corresponding extended CMV matrices by $\mathcal{E}$ and define $\mathcal{L}$, $\mathcal{M}$ as
\begin{equation*}
\mathcal{L}=\bigoplus_{j\in \mathbb{Z}}\Theta_{2j}, \mathcal{M}=\bigoplus_{j\in \mathbb{Z}}\Theta_{2j+1}.
\end{equation*}

The analogous factorization of $\mathcal{E}$ is given by $\mathcal{E}=\mathcal{L}\mathcal{M}$.

\medskip

Let $\mathcal{E}_{[a,b]}$ denote the restriction of an extended CMV matrix to the finite interval [a,b], defined by
 $$\mathcal{E}_{[a,b]}=(P_{[a,b]})^{\ast}\mathcal{E}P_{[a,b]},$$
where $P_{[a,b]}$ is the projection $\ell^{2}(\mathbb{Z})\rightarrow \ell^{2}([a,b])$. $\mathcal{L}_{[a,b]}$ and $\mathcal{M} _{[a,b]}$ are defined similarly.

With $\beta,\gamma \in \partial \mathbb{D}$, define the sequence of Verblunsky coefficients
\begin{equation*}
\tilde{\alpha_n}=
\begin{cases}
\beta, \qquad & n =a;\\
\gamma,& n =b;\\
\alpha_n,& n \notin  \{a,b\}.
\end{cases}
\end{equation*}

The corresponding operator is defined by $\tilde{\mathcal{E}}$ and we define
$\mathcal{E}^{\beta,\gamma}_{[a,b]}=P_{[a,b]}\tilde{\mathcal{E}}\left(P_{[a,b]}\right)^*$.
$\mathcal{E}^{\beta,\gamma}_{[a,b]}$ is unitary whenever $\beta,\gamma \in \partial\mathbb{D}$, which can be verified.
Then, for $z\in\mathbb{C}$, $\beta,\gamma\in\partial\mathbb{D}$, the polynomials are defined by
$$\Phi_{[a,b]}^{\beta,\gamma}(z):=\det{\left(z-\mathcal{E}_{[a,b]}^{\beta,\gamma}\right)}, \qquad
\phi_{[a,b]}^{\beta,\gamma}(z):=(\rho_a\cdots\rho_b)^{-1}\Phi_{[a,b]}^{\beta,\gamma}(z).$$

Since the equation $\mathcal{E}\psi=z \psi$ is equivalent to $(z\mathcal{L}^{\ast}-\mathcal{M})\psi=0$. Then the associated finite-volume Green's functions are as follows
\begin{equation*}
G^{\beta,\gamma}_{\omega,[a,b]}(z)=(z[\mathcal{L}^{\beta,\gamma}_{\omega,[a,b]}]^{*}-[\mathcal{M}^{\beta,\gamma}_{\omega,[a,b]}])^{-1},
\end{equation*}
and
\begin{equation*}
G^{\beta,\gamma}_{\omega,[a,b]}(j,k,z)=\langle\delta_{j},G^{\beta,\gamma}_{\omega,[a,b]}(z)\delta_{k}\rangle,~~j,k\in [a,b].
\end{equation*}
By \cite[Proposition 3.8]{Kr}, the Green's function has the following expression
\begin{equation}{\label{2.16}}
\left|G_{\omega,[a,b]}^{\beta,\gamma}(j,k;z)\right|=\frac{1}{\rho_j\rho_k}
\left|\frac{\phi_{\omega,[a,j-1]}^{\beta,\alpha_{j-1}}(z)\phi_{\omega,[k+1,b]}^{\alpha_{k+1},\gamma}(z)}
{\phi_{\omega,[a,b]}^{\beta,\gamma}(z)}\right|.
\end{equation}

Now one can get the following Green's function estimate.
\begin{lemma}\label{e.greenfunction01} (\cite[Lemma 3.5]{Wang01})
Assume that for $n$ large enough and any $\varepsilon>0$, the following inequality holds,
$$\frac{1}{n}\log{\|M_n^z(x)\|}\geq L_n(z)-\varepsilon.$$
Then we have that, for any $\beta_0,\gamma_0\in\partial\mathbb{D}$, there exist $\beta\in\{\beta_0,-\beta_0\}$ and $\gamma\in\{\gamma_0,-\gamma_0\}$,
such that
\begin{equation}{\label{2.17}}
\left|G_{\omega,[0,n)}^{\beta,\gamma}(j,k;z)\right|\leq e^{-|j-k|L_n(z)+C\varepsilon n},
\end{equation}
for all $j,k\in[0,n)$, $z\in\partial\mathbb{D}\backslash \sigma\left(\mathcal{E}_{\omega,[0,n)}^{\beta,\gamma}\right)$.
\end{lemma}

\subsection{Elimination of double resonances and semi-algebraic sets}

The following lemmas are basically contained in \cite[Sect. 3]{Bourgain03}. For convenience, we only state the two lemmas here.
\begin{lemma}\label{3.3}(\cite[Lemma 3.3]{Bourgain03})
Let $T_{\omega}:\mathbb{T}^{2}\rightarrow \mathbb{T}^{2}$ be the $\omega$-skew-shift, and  $S\subset \mathbb{T}^{4}\times \mathbb{C}$ be a semi-algebraic set of degree at most $B$ such that
\begin{equation}\label{3.5}
\mathrm{mes}(\mathrm{Proj}_{\mathbb{T}^{4}}S)<e^{-B^\sigma}~~ for~some ~~\sigma>0.
\end{equation}
Suppose $M$ and $B$ are related by the inequalities
\begin{equation}
\log\log M\ll \log B \ll \log M,
\end{equation}
then,
\begin{equation}\label{3.9}
\mathrm{mes}\left[(y_{0},\omega)\in \mathbb{T}^2|(y_{0},\omega, T^{j}_{\omega}(0,y_{0}))\in \mathrm{Proj}_{\mathbb{T}^{4}}S ~~ \mathrm{for~some} ~~j\sim M\right]<M^{-10^{-8}}.
\end{equation}
\end{lemma}
\begin{lemma}\label{3.6}(\cite[Lemma 3.4]{Bourgain03})
Let $T_{\omega}$ be the $\omega$-skew-shift, $\omega$ satisfying
\begin{equation}\label{3.40}
\|k\omega\|\geq c_{\varepsilon}|k|^{-1-\varepsilon}~for~all~k\in \mathbb{Z}, 0<|k|<N.
\end{equation}
Then, denoting
$$u_{N_{0}}(x,y)=\frac{1}{N_{0}}\log \left\|(\sqrt{z})^{-n} \prod^{0}_{j=N_{0}-1}\frac{1}{\rho(T^{j}_{\omega}(x,y))}
\left(
\begin{matrix}
\sqrt{z} &\frac{ -\overline{\alpha}(T^{j}_{\omega}(x,y))}{\sqrt{z}}\\
-\alpha(T^{j}_{\omega}(x,y))\sqrt{z} & \frac{1}{\sqrt{z}}
\end{matrix}
\right)\right\|,$$
there exist constants $\sigma>0$ and $C>1$ so that for $N>N^{C}_{0}$ one has the uniform bound
\begin{equation}\label{3.41}
\left\|\frac{1}{N}\sum^{N-1}_{j=0}\left|u_{N_0}\circ T^{j}_{\omega}-\int_{\mathbb{T}^2}u_{N_0}(x,y)dxdy \right| \right\|_{L^{\infty}(\mathbb{T}^2)}<N^{-\sigma}_{0}.
\end{equation}
\end{lemma}

According to the previous lemmas, the uniform upper bound on the norm of the monodromy matrices can be obtained, as follows
\begin{coro}(\cite[Corollary 3.5]{Bourgain03})\label{3.7}
Assume $\omega$ satisfies \eqref{3.40}. For any $N>N^{C}_0$, there is a uniform estimate for all $z\in \mathbb{C}$,
\begin{equation}\label{3.54}
\sup_{(x,y)\in\mathbb{T}^2}\frac{1}{N}\log\|M_{N}(x,y,z)\|<L_{N_0}(z)+N^{-\sigma}_0.
\end{equation}
\end{coro}

In \cite{Bourgain01}, Section 11 the elimination of double resonances at a fixed point is proved based on the following fact for the restriction $H_{[0,n-1]}$ of a Schr\"{o}dinger operator $H$, since it is self-adjoint, then
$$\mathrm{dist}(E,\mathrm{spec}(H_{[0,n-1]}))=\|(H_{[0,n-1]}-E)^{-1}\|^{-1}~~\textrm{for}~~ E\notin \mathrm{spec}(H_{[0,n-1]}).$$

Although $\mathcal{E}_{[0,n-1]}$ is not self-adjoint, even not normal, we still can get the same estimate as in Schr\"{o}dinger case by using the result from \cite{Davis-Simon}.

\begin{lemma}(\cite[Theorem 1]{Davis-Simon})\label{th.distance}
Let $A_{n}$ be the set of pairs $(A,z)$, where $A$ is an $n\times n$ matrix, $z\in \mathbb{C}$ with
$$|z|\geq \|A\|$$
and
$$ z \notin \mathrm{spec}(A).$$
Then
$$\sup_{A_{n}}\mathrm{dist}(z,\mathrm{spec}(A))\|(z-A)^{-1}\|= \mathrm{cot}(\frac{\pi}{4n}).$$
\end{lemma}

Fix $\varepsilon>0$ small and let $\omega \in \Omega_\varepsilon$, see \eqref{e.diophantine}. By Theorem \ref{mainth1}, we have
\begin{equation}\label{e.lyapunov01}
\inf_{z}L(z)>c_{0}>0,
\end{equation}
where $c_0$ is an absolute constant.

\begin{lemma}\label{le.double01}
Fix $y_0 \in \mathbb{T}$, $\beta, \gamma \in \partial \mathbb{D}$, and  $\varepsilon>0$ small. Let constant $C_{1}\geq1$ and $N$ be an arbitrary positive integer. Define $S=S_{N}\subset \mathbb{T}^{4}\times \partial \mathbb{D}$ to be the set of those $(\omega,y_{0},x,y,z)$ for which there exists some $N_{1}<N^{C_1}$ so that
\begin{equation}\label{3.56}
\|k\omega\|\geq\varepsilon|k|^{-1}(1+\log k)^{-2}~ \textrm{ for all }~ k\in \mathbb{Z}, 0<k<N,
\end{equation}
\begin{equation}\label{3.57}
\left\|(\mathcal{E}^{\beta,\gamma}_{[-N_{1},N_1]}(\omega,0,y_0)-z)^{-1} \right\|>e^{C_{2}N},
\end{equation}
\begin{equation}\label{3.58}
\frac{1}{N}\log\|M_{N}(\omega,x,y,z)\|<L_{N}(\omega,z)-\frac{c_0}{10}.
\end{equation}
Here $c_0$ is the constant from \eqref{e.lyapunov01} and $C_2$ is a sufficiently large constant depending on $\alpha$. Then
\begin{equation}\label{3.59}
\mathrm{mes}(\mathrm{Proj}_{\mathbb{T}^{4}}S)\lesssim e^{-\frac{1}{2}N^{\sigma}}.
\end{equation}

Moreover, $S$ is contained in a set $S^{'}$ satisfying the measure estimate \eqref{3.59} and which is semi-algebraic of degree at most $N^C$ for some constant $C$ depending on $\alpha, \varepsilon,C_1, C_2$.
\end{lemma}
\Proof Fix $N$ large, denote
\begin{equation*}
\Lambda_{\omega}^{\beta,\gamma}=\bigcup_{N_1<N^{C_2}}\mathrm{spec}(\mathcal{E}^{\beta,\gamma}_{[-N_1,N_1]}(\omega,0,y_0))
\end{equation*}
and
\begin{equation*}
\Upsilon_{\omega}^{\beta,\gamma}=\{(x,y)\in\mathbb{T}^2: \max_{z_1 \in \Lambda_{\omega}^{\beta,\gamma}}|\frac{1}{N}\log\|M_N(x,y;z_1)\|-L_{N}(z_1)|>\frac{c_0}{20}\},
\end{equation*}
By the large deviation estimate, we have
\begin{equation}\label{e.double01}
\mathrm{mes}\Upsilon_{\omega}^{\beta,\gamma}\lesssim e^{-N^\sigma},
\end{equation}
then \eqref{3.59} will follow from the fact that
\begin{equation}\label{e.double02}
\mathrm{Proj}_{\mathbb{T}^{4}}S \subset\{(y_0,\omega,x,y)\in \mathbb{T}^4: \|k\omega\|\geq\varepsilon|k|^{-1}(1+\log k)^{-2}; (x,y)\in \Upsilon_{\omega}^{\beta,\gamma}\}.
\end{equation}

Recall the large deviation estimate for $n=N$ holds under the condition \eqref{3.56} on $\omega$. Assume \eqref{3.57} and \eqref{3.58} holds for some $z\in\partial\mathbb{D}$.
Due to Lemma \ref{th.distance}, we have
\begin{equation*}
\mathrm{dist}(z,\mathrm{spec}(\mathcal{E}^{\beta,\gamma}_{[-N_1,N_1]}(\omega,0,y_0))) \|(z-\mathcal{E}^{\beta,\gamma}_{[-N_1,N_1]}(\omega,0,y_0))^{-1}\| \le c(N),
\end{equation*}
which means that
\begin{equation*}
\mathrm{dist}(z,\mathrm{spec}(\mathcal{E}^{\beta,\gamma}_{[-N_1,N_1]}(\omega,0,y_0))) \le c(N) \|(z-\mathcal{E}^{\beta,\gamma}_{[-N_1,N_1]}(\omega,0,y_0))^{-1}\|^{-1}.
\end{equation*}
Thus, there exists $z_{1}\in \mathrm{spec}(\mathcal{E}^{\beta,\gamma}_{[-N_1,N_1]}(\omega,0,y_0))$ which do not depend on $(x,y)$, such that
\begin{equation}\label{e.double03}
|z-z_1|<e^{-C_2N}.
\end{equation}
It follows from \eqref{e.double03} with sufficiently large $C_2$ and \eqref{3.58} that
\begin{equation}\label{e.double04}
|\frac{1}{N}\log\|M_N(x,y;z_1)\|-L_{N}(z_1)|>\frac{c_0}{20},
\end{equation}
the measure of the set of $(x,y)\in \mathbb{T}^2$ for which \eqref{e.double04} holds for fixed $z_1$ does not exceed $e^{-N^\sigma}$. This proves that
\begin{equation*}
\mathrm{mes}(\mathrm{Proj}_{\mathbb{T}^{4}}S)\leq 4N_{1}^{2} e^{-N^\sigma}\lesssim e^{-\frac{1}{2}N^\sigma}.
\end{equation*}

\medskip

Furthermore, we still need to prove that the conditions \eqref{3.57} and \eqref{3.58} can be replaced by inequalities involving only polynomials of degree at most $N^C$, without increasing the measure estimate \eqref{3.59} by more than a factor or two. The full details are not described here, see \cite{Bourgain03,Bourgain01}.
\qedbox

\subsection{The proof of Theorem \ref{mainth2} }

Let $\omega\in \Omega_\varepsilon$, see \eqref{e.diophantine}. For large $N$, let $S_N$ be as in Lemma \ref{le.double01}. Then Lemma \ref{3.3} applies to $S_N$ and setting $\overline{N}=e^{(\log N)^{2}}$ it follows that
\begin{equation}\label{3.63}
\mathrm{mes}\left [(y_{0},\omega)\in \mathbb{T}^2 \mid(y_{0},\omega, T^{j}_{\omega}(0,y_{0}))\in \mathrm{Proj}_{\mathbb{T}^{4}}(S_N) ~\mathrm{ for~ some} ~j\sim \overline{N}\right]<\overline{N}^{-10^{-8}}.
\end{equation}

Let $\mathcal{B}_N$ denote the set on the left-hand side of \eqref{3.63} and define
\begin{equation*}
\mathcal{B}^{(0)}:=\lim\sup_{N\rightarrow \infty}\mathcal{B}_N.
\end{equation*}
Thus $\mathrm{mes}(\mathcal{B}^{(0)})=0$. Since $T^{\ell}(x,y)=(x,0)+T^{\ell}(0,y)(\mathrm{mod} 1)$, this construction applied to the potential $\alpha(x+\cdot,\cdot)$ instead of $\alpha$ produces a set $\mathcal{B}^{(x)}$ of measure zero. Finally, set
\begin{equation*}
\mathcal{B}:=\left\{(\omega,x,y)|(y,\omega)\in \mathcal{B}^{(x)} \right\},
\end{equation*}
which is again of measure zero. It is for all $(\omega,x,y)\in (\Omega_{\varepsilon}\times \mathbb{T}^2)\backslash \mathcal{B}$ that we shall prove localization.

\medskip

Fix such a choice of $(\omega,x,y)$ and any $z\in \partial \mathbb{D}$. $\xi=\{\xi_n\}_{n\in \mathbb{Z}}$ satisfies the equation
\begin{equation*}
\mathcal{E}_{\omega}\xi=z\xi
\end{equation*}
where
\begin{equation}\label{3.64}
|\xi_n|\lesssim |n|^C ~\mathrm{for~ all}~n \in \mathbb{Z}.
\end{equation}
Furthermore, we normalize $|\xi_0|=|\xi_1|=1$.

By our choice of $(\omega,x,y)$, fix some large integer $N$ and if there is $N_{1}<N^{C_1}$ such that
\begin{equation}\label{3.57*}
\left\|(\mathcal{E}^{\beta,\gamma}_{[-N_{1},N_1]}-z)^{-1} \right\|>e^{C_{2}N}
\end{equation}
holds for $\beta,\gamma \in \partial\mathbb{D}$,
according to Lemma \ref{le.double01}, then \eqref{3.58} must fail, which implies
\begin{equation*}
\frac{1}{N'}\log\|M_{N'}(T^{j}_{\omega}(x,y);z)\|>L(z)-\frac{c_0}{10}
\end{equation*}
for all $N'\sim N$ and $j\sim \overline{N}=e^{(\log N)^2}$. It follows from the avalanche principle that if $\frac{\overline{N}}{2}<|j|<\overline{N}$
 and $N^{2}<N_{2}<\frac{\overline{N}}{10}$, then also
\begin{equation}\label{3.65}
\frac{1}{N_2}\log\|M_{N_2}(T^{j}_{\omega}(x,y);z)\|>L(z)-\frac{c_0}{10}.
\end{equation}

\medskip

As usual, let
\begin{equation*}
G^{\beta,\gamma}_{\Lambda}(\omega,x,y;z):=(\mathcal{E}^{\beta,\gamma}_{\Lambda}(\omega,x,y)-z)^{-1}
\end{equation*}
be the Green's function. $\mathcal{E}_{\omega,\Lambda}$ denotes the restriction of $\mathcal{E}_{\omega}$ to the interval $\Lambda$. Consider
 intervals
 $$\Lambda=\left[j,j+\frac{\overline{N}}{10}\right],$$
where $\frac{\overline{N}}{2}<|j|<\overline{N}$.

By definition of $G_{\Lambda}$ and because of \eqref{3.64}, it will suffice to prove that for fixed constant $c_{1}>0$,
\begin{equation}\label{3.66}
\max_{\ell\in\partial\Lambda}\left|G^{\beta,\gamma}_{\Lambda}(\omega,x,y;z)(k,\ell)\right|\lesssim \mathrm{exp}(-c_{1}\overline{N})~\mathrm{for~all}~
k\in\Lambda,
\end{equation}
with $\mathrm{dist}(k,\partial\Lambda)>\frac{1}{4}|\Lambda|$. Once \eqref{3.66} holds, then one can get the required exponential decay property as follows,
\begin{equation*}
\begin{split}
|\xi_{k}| <&|G^{\beta,\gamma}_{\Lambda}(k,\ell)||\xi_{\ell}|+|G^{\beta,\gamma}_{\Lambda}(k,\ell+\frac{\overline{N}}{10})||\xi_{\ell+\frac{\overline{N}}{10}}|\\
\lesssim & \overline{N}^C\mathrm{exp}(-c_{1}\overline{N})\\
\lesssim & \exp(-\frac{c_{1}\overline{N}}{2}),
\end{split}
\end{equation*}
where we used \eqref{C1} in the first step and \eqref{3.64} in the second step.

\medskip

According to the above statement, we first need to prove that for $N_{1}<N^{C_1}$
$$\|G^{\beta,\gamma}_{[-N_{1},N_{1}]}\|>e^{C_{2}N}.$$
Recalling $|\xi_{0}|=1$, by \cite[Lemma 3.9]{Kr}, we have
\begin{equation}
1 \leqslant\|\xi_{[-N_{1},N_{1}]}\| \leqslant 2\|G^{\beta,\gamma}_{[-N_{1},N_{1}]}\|(|\xi_{-N_1}|+|\xi_{-N_{1}+1}|+|\xi_{N_{1}-1}|+|\xi_{N_1}|).
\end{equation}
It suffices to show that
\begin{equation*}
|\xi_{-N_1}|+|\xi_{-N_{1}+1}|+|\xi_{N_{1}-1}|+|\xi_{N_1}|\le \frac{1}{2}e^{-C_{2}N}
\end{equation*}
for some $N_{1}\thicksim N^{C_1}$.

Now, we give the estimation of $|\xi_{N_{1}}|$. According to Lemma \ref{3.6}, there exists a $j\sim N^{C_1}$ such that
\begin{equation}\label{e.LDT0}
\left| u_{4C_2 N}\circ T^{j}_{\omega}-L_{4C_2 N}\right|\lesssim N^{-\sigma},
\end{equation}
where $N$ is large enough such that $N^{-\sigma}<\frac{c_0}{10}$ with $c_{0}>0$
is the lower bound of $L(z)$. Integrating \eqref{e.LDT0} over $\T^2$, we have
\begin{equation*}
\frac{1}{4C_2N}\log\|M_{4C_2 N}(T^{j}_{\omega}(x,y);z)\|>L(z)-\frac{c_0}{10}.
\end{equation*}
By Lemma \ref{e.greenfunction01}, it follows that \eqref{3.66} holds. Combining with \eqref{C1} and \eqref{3.64}, one can obtain that
$$
|\xi_{N_{1}}| \le \frac{1}{8} e^{-C_2 N}.
$$
Obviously, we also have $
|\xi_{N_{1}-1}|, |\xi_{-N_{1}}|, |\xi_{-N_{1}+1}| \le \frac{1}{8} e^{-C_2 N}.
$

Thus, the assertion holds.
\qedbox

\vskip1cm

\begin{appendix}

\section{Avalanche principle}
The tool to exploit the structure of $M_{n}^{z}$ more carefully is the ``avalanche principle" from Goldstein and  Schlag \cite{Goldstein01}.
\begin{lemma}\label{le.AP}
Let $A_{1}, \cdots, A_n$ be a sequence of arbitrary unimodular 2$\times $2-matrices. Suppose that
\begin{equation}\label{2.8}
\min_{1\leq j\leq n}\|A_j\|\geq \mu \geq n~ and
\end{equation}
\begin{equation}\label{2.9}
\max_{1\leq j\leq n}\left[\log\|A_{j+1}\|+\log\|A_j\|-\log\|A_{j+1}A_j\|\right]\leq\frac{1}{2}\log \mu.
\end{equation}
Then
\begin{equation}\label{2.10}
\left| \log\|A_{n}\cdot \ldots \cdot A_1\|+\sum^{n-1}_{j=2}\log\|A_j\|-\sum^{n-1}_{j=1}\log\|A_{j+1}A_j\| \right|<C\frac{n}{\mu}.
\end{equation}
\end{lemma}
\section{Subharmonic functions}
Denote
$$\mathcal{D}_{h}:=\{t\in\mathbb{C}:1-h<|t|<1+h\},$$
suppose $u$ is subharmonic on $\mathcal{D}_{h}$, with $\sup_{\mathcal{D}_{h}}|u|\leq N$. By Riesz's representation theorem, there is a positive measure $\mu$ with $\mathrm{supp}(\mu)\subset \mathcal{D}_{\frac{h}{2}} $ and a harmonic function $H$ such that for any $t \in \mathcal{D}_{\frac{h}{2}}$
\begin{equation}\label{e.subharmonic01}
u(t)=\int \log|t-\zeta|d\mu(\zeta)+H(t),
\end{equation}
where
\begin{equation*}
\mu(\mathcal{D}_{\frac{h}{2}})+\|Y\|_{L^{\infty}(\mathcal{D}_{\frac{h}{4}})}\leq C_{h}N.
\end{equation*}
\begin{lemma}\label{le.subharmonic01} (\cite[Lemma 2.5]{Bourgain03})
Let $u: \mathbb{T}^2\rightarrow \mathbb{R}$ satisfy $\|u\|_{L^{\infty}(\mathbb{T}^2)}\leq1$. Assume that $u$ extends as a separately subharmonic function in each variable to a neighborhoods of $\mathbb{T}^2$ such that for some $N\geq1$ and $\rho>0$,
$$\sup_{(w_1,w_2)\in \mathcal{D}_{h_1}\times\mathcal{D}_{h_2}}|u(w_1,w_2)|\leq N.$$
Furthermore, suppose that $u=u_{0}+u_{1}$ on $\mathbb{T}^2$,
$$\|u_{0}-\langle u\rangle\|_{L^{\infty}(\mathbb{T}^2)} \leq\varepsilon_0, ~ \|u_1\|_{L^{1}(\mathbb{T}^2)}\leq\varepsilon_1$$
with $0<\varepsilon_0,\varepsilon_1<1$, and $\langle u\rangle:=\int_{\mathbb{T}^2}u(x,y)dxdy$. Then for any $\delta>0$, and
\begin{equation*}
B=\varepsilon_{0}\log(\frac{N}{\varepsilon_1})+N^{\frac{3}{2}}\varepsilon_{1}^{\frac{1}{4}},
\end{equation*}
one can get
\begin{equation}
\mathrm{mes}\left[(x,y)\in\mathbb{T}^2: |u(x,y)-\langle u\rangle|>B^{\delta}\log(\frac{N}{\varepsilon_1})\right]\leq C_{\rho}N^{2}\varepsilon_1^{-1}\exp(-c_{\rho}B^{\delta-\frac{1}{2}}).
\end{equation}
\end{lemma}
\begin{lemma}\label{le.subharmonic02}(\cite[Lemma 2.6]{Bourgain03})
Let $u$ satisfy the conditions in Lemma \ref{le.subharmonic01}, so that for some $\rho>0$,
$$\sup_{(w_1,w_2)\in \mathcal{D}_{h_1}\times\mathcal{D}_{h_2}}|u(w_1,w_2)|\leq 1.$$
Fix $\varepsilon>0$ small, $\omega\in \Omega_{\omega}$, for any $\delta>0$, there exist constants $c, C$ (only depending on $\rho$, $\delta$, $\varepsilon$),
\begin{equation}
\mathrm{mes}\left[(x,y)\in\mathbb{T}^2: |\frac{1}{K}\sum_{k=1}^{K}u\circ T^{k}_{\omega}(x,y)-\langle u\rangle|>K^{2\delta-\frac{1}{12}}\right]\leq C\exp(-cK^{\delta}).
\end{equation}
where $K$ is any positive integer.
\end{lemma}

\section{Restriction of eigenequations}

In the Schr\"{o}dinger case, by restriction the eigenequation $(H-E)\xi=0$ to a finite interval $\Lambda=[a,b]$, one can get two boundary terms, which yields the identity $\xi(n)=-G_{\Lambda}^{E}(n,a)\xi(a-1)-G_{\Lambda}^{E}(n,b)\xi(b+1)$.

But in the CMV case, the analog of this formula depends on the parity of the endpoints of the finite interval. Concretely, if $\psi$ is a solution of difference equation $\mathcal{E}\psi=z\psi$,
we define
\begin{equation*}
\tilde{\psi}(a)=
\begin{cases}
(z\overline{\beta}-\alpha_a)\psi(a)-\rho_a\psi(a+1), \qquad & a\text{ is even,}\\
(z\alpha_a-\beta)\psi(a)+z\rho_a\psi(a+1),& a \text { is odd,}
\end{cases}
\end{equation*}
and
\begin{equation*}
\tilde{\psi}(b)=
\begin{cases}
(z\overline{\gamma}-\alpha_b)\psi(b)-\rho_b\psi(b-1), \qquad & b\text{ is even,}\\
(z\alpha_b-\gamma)\psi(b)+z\rho_{b-1}\psi(b-1),& b \text { is odd}.
\end{cases}
\end{equation*}
Then for $a<n<b$, we have
\begin{equation}{\label{C1}}
\psi(n)=G_{[a,b]}^{\beta,\gamma}(n,a;z)\tilde{\psi}(a)+G_{{[a,b]}}^{\beta,\gamma}(n,b;z)\tilde{\psi}(b).
\end{equation}
\end{appendix}

 \vskip1cm

\noindent{$\mathbf{Acknowledgments}$}

This work was supported by the NSFC (No. 11571327).

\end{document}